# Fuzzy Finite Element Solution of Uncertain Convection-Diffusion Heat Transfer for A Rectangular Plate

Sudipta Priyadarshini[1] and Sukanta Nayak[2], Paresh Kumar Panigrahi[3]

[1, 2, 3]Department of Mathematics, School of Advanced Sciences, VIT-AP University, Amaravati, Andhra Pradesh

## ABSTRACT

Convection-diffusion of heat transfer is one of the important phenomena in fluid flow and industrial problems. The involved parameters, boundary conditions, and material properties are greatly affecting the same. As such, the uncertainness of these parameters, conditions, and properties cannot be ignored. Therefore, the present research work focuses numerical approach to study the uncertain convection-diffusion of heat transfer problem. Here, the finite element method is employed with fuzzy uncertainties to investigate the uncertain temperatures for a plate problem. The uncertainty is considered with a 10% of error and then corresponding fuzzy numbers are generated. These fuzzy numbers are used in the governing differential equation and boundary conditions to get the nodal temperatures. A different combination of fuzzy parameters is considered and the obtained results are reported. Furthermore, the sensitiveness of these parameters is reported in a case study.

## 1. INTRODUCTION

Uncertainty plays an important role in various fields of science and technology. These uncertainties are arising due to experimental error, vague information, and insufficient knowledge of objects involved in the system. In convection-diffusion heat transfer problems the uncertainties are encountered due to the heat constants, coefficients and boundary conditions. To quantify the uncertainness of the problem first we need to understand the system. As such, first few important research works are reported in ( [1], [2], [3], [4], [5]) to understand the system. Then to handle the uncertainty and investigate the uncertain system researchers adopted probabilistic approach and quantified the temperatures. In this regards, Deng and Liu [6] considered Monte Carlo method to solve a bio heat transfer problem. Wu [7] used a Monte Carlo approach to simulate transient radiative transfer in a refractive planar media exposed to collimated pulse irradiation. Also, many of the articles are present that uses Monte Carlo method for the same. These approaches need large data set to analyse the thermal problems. Hence, we need an alternate approach to handle small data sets as well as method to solve the same.

In view of the above, numerical techniques can be employed to investigate the thermal problems. Well known numerical techniques viz. Finite Difference Method (FDM), Finite Volume Method (FVM), and Finite Element Method (FEM) can be used to quantify the temperatures. Due to the better usability, FEM is more familiar for the said systems. Couple of important researches are included here to know the importance of FEM. Edward and Robert [8] have used the Finite Element Method to investigate heat conduction problem. Further, Onate et al. [9] have used Galerkin finite element method for solving convection-diffusion problems with sharp gradients with calculus. But, from the above literature review it is seen the problems are discussed without uncertainties as well as the large data sets are taken that may be cost worthy to handle. In this context, fuzzy sets can be implemented that can handle small data set. Zadeh [10] introduced fuzzy sets and fuzzy number to compute fuzzy uncertainties involved in the system. Dong and Wong [11] have discussed the extension principle's computational aspect by using the fuzzy weighted average method (FWAM). Chakraverty and Nayak [12] have described interval/fuzzy uncertainty for solving conduction problems. The interval arithmetic is extended to the fuzzy number and the fuzzy finite element method is developed. Then, Nayak and Chakraverty [13] took the modified Fuzzy Finite Element Method (FFEM) for solving the conduction-convection system under a square plate. Furthermore, Nayak et al. [14] investigated the temperature distribution in a conduction-convection system under ambiguous environments. Whereas, Wang and Qiu [15] take a modified version of the intervals, that is an Interval Parameter Perturbation Method (IPPM) and a sensitivity-based interval analysis method which is also used for solving heat conduction problems. Several numerical strategies were also presented in Li et al. [16] to solve the convective diffusion equation with random diffusivity and periodic boundary conditions. Recently, Nayak [17] has taken the Interval Finite Element Method (IFEM) for solving parametric uncertainties in convection-diffusion problems. The above literatures provide the importance of the role of uncertainty with the solution strategies. So, taking this into account, this paper presents a convection-diffusion heat flow on a rectangular pate with uncertain involved parameters. The same is investigated with fuzzy finite element method. Then the obtained results are analysed. Based on the variation of the uncertain temperatures the involved parameters are studied.

---

[1] Email: sudipta.grc@gmail.com (Sudipta Priyadarshini)
[2] Email: sukantgacr@gmail.com (Sukanta Nayak)
3 Email: jrpareshkp@gmail.com (Paresh Kumar Panigrahi)

## 2. Fuzzy finite element formulation of the convection-diffusion problem

Consider the following three-dimensional governing transient convection-diffusion equation

$$\frac{\partial \varphi}{\partial t} + \sum_{i=1}^{3} \frac{\partial \varphi}{\partial x_i} + \varphi \sum_{i=1}^{3} \frac{\partial u_i}{\partial t} - \sum_{i=1}^{3} \frac{\partial}{\partial x_i}\left(k \frac{\partial \varphi}{\partial x_i}\right) + Q = 0 \tag{1}$$

where, $\varphi$ is the temperature, $k$ is the diffusion coefficient, $u_i$ are the convection velocity components and Q is the source. In Eq. (1), the first term is transient, the second and third terms are convection terms and the fourth term represents the diffusion. The Galerkin weak form of the differential equation (1) is

$$\int_{\Omega} w \cdot \frac{\partial \varphi}{\partial t} d\Omega + \int_{\Omega} w \cdot u_i \frac{\partial \varphi}{\partial x_i} d\Omega + \int_{\Omega} w \cdot \varphi \frac{\partial u_i}{\partial x_i} d\Omega - \int_{\Omega} w \cdot \frac{\partial}{\partial x_i}\left(k \frac{\partial \varphi}{\partial x_i}\right) d\Omega + \int_{\Omega} w \cdot Q \; d\Omega = 0 \tag{2}$$

where, $w$ is defined as the weight function.

Here, the shape functions nothing but the weight functions. Eq. (2) can be rewritten as

$$\int_{\Omega} w \cdot \nabla \varphi \; d\Omega + \int_{\Omega} w \cdot u_i \nabla \varphi \; d\Omega + \int_{\Omega} w \cdot \varphi \nabla u \; d\Omega - \int_{\Omega} w \cdot \nabla(k \nabla \varphi) d\Omega + \int_{\Omega} w \cdot Q \; d\Omega = 0. \tag{3}$$

If one dimension case is considered then Eq. (1) can be transferred to

$$\frac{\partial \varphi}{\partial t} + u_1 \frac{\partial \varphi}{\partial x_1} + \varphi \frac{\partial u_i}{\partial x_1} - \frac{\partial}{\partial x_1}\left(k \frac{\partial \varphi}{\partial x_1}\right) + Q = 0. \tag{4}$$

Further, for constant convection velocity $u_1$, Eq. (4) can be represented as

$$\frac{\partial \varphi}{\partial t} + u_1 \frac{\partial \varphi}{\partial x_1} - \frac{\partial}{\partial x_1}\left(k \frac{\partial \varphi}{\partial x_1}\right) + Q = 0. \tag{5}$$

A one-dimensional convection equation without a source term may be obtained by neglecting the diffusion and source terms as follows

$$\frac{\partial \varphi}{\partial t} + u_1 \frac{\partial \varphi}{\partial x_1} = 0. \tag{6}$$

After the weak formulation of the governing differential equation, fuzzy parameters can be implemented to quantify the uncertain field variables.

### 2.1 Fuzzy Finite Element Method (FFEM)

The idea given in (Chakraverty and Nayak [12]) is used here to solve an uncertain convection-diffusion heat transfer problem for a rectangular plate. For the sake of completeness, few definitions are included below.

**Definition 1**

A fuzzy number is defined as a piecewise continuous convex and normalized fuzzy set. It can be denoted as $A` = [\,a^L, a^M, a^R]$. Here $A`$ is a Triangular Fuzzy Number (TFN). The membership function of $A`$ is

$$\mu_{A`}(\xi) = \begin{cases} 0, \text{otherwise} \\ \mu_L, a^L \leq \xi \leq a^M \\ \mu_R, a^M \leq \xi \leq a^R \end{cases} \qquad (7)$$

where, $\mu_L = \frac{\xi - a^L}{a^M - a^L}$ and $\mu_R = \frac{a^R - \xi}{a^R - a^N}$.

The TFN $A` = [a^L, a^M, a^R]$ may be transformed into interval form by using $\alpha - $ cut as follow

$$A` = [a^L, a^M, a^R] = [a^L + (a^M - a^L)\alpha, a^R - (a^M - a^R)\alpha], \quad \alpha \in [0,1]. \qquad (8)$$

The TFN $A`$ is shown in the Fig. 1.

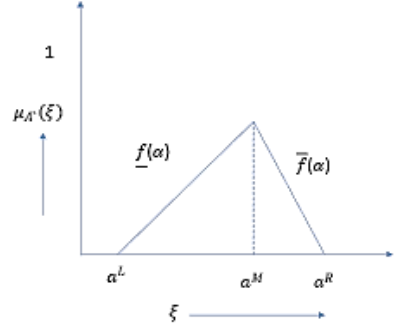

Fig. 1. Triangular Fuzzy Number (TFN)

**Definition 2**

Consider two fuzzy numbers $\xi` = [\overline{\xi}(\alpha), \underline{\xi}(\alpha)]$ and $\gamma` = [\overline{\gamma}(\alpha), \underline{\gamma}(\alpha)]$ and a scaler $\ell$ then

(a) $\xi` = \gamma`$ if and only if $\underline{\xi}(\alpha) = \underline{\gamma}(\alpha)$ and $\overline{\xi}(\alpha) = \overline{\gamma}(\alpha)$;

(b) $\xi` + \gamma` = [\underline{\xi}(\alpha) + \underline{\gamma}(\alpha), \overline{\xi}(\alpha) + \overline{\gamma}(\alpha)]$;

(c) $\ell\xi` = \begin{cases} [\ell\underline{\xi}(\alpha), \ell\overline{\xi}(\alpha)], \ell \geq 0 \\ [\ell\overline{\xi}(\alpha), \ell\underline{\xi}(\alpha)], \ell < 0 \end{cases}$.

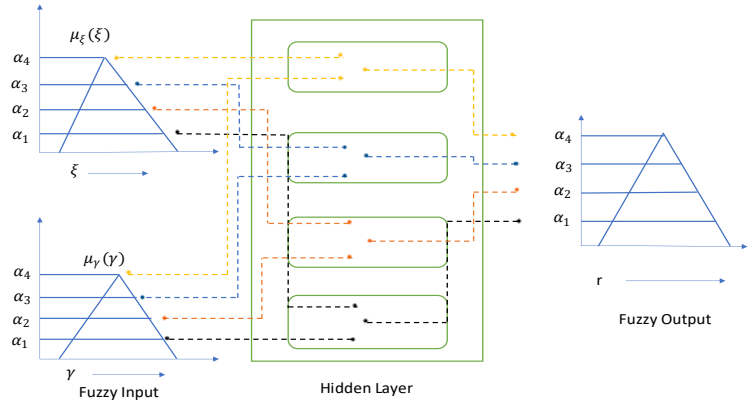

Fig.2 Model diagram for fuzzy Finite Element Method

**Definition 3**

The fuzzy arithmetic is also redefined ([13]) as

(I) $[\overline{\xi}(\alpha), \underline{\xi}(\alpha)] + [\overline{\gamma}(\alpha), \underline{\gamma}(\alpha)] = [\min\{\lim_{n\to\infty} m_1 + \lim_{n\to\infty} m_2, \lim_{n\to 1} m_1 + \lim_{n\to 1} m_2\}, \max\{\lim_{n\to\infty} m_1 + \lim_{n\to\infty} m_2, \lim_{n\to 1} m_1 + \lim_{n\to 1} m_2\}]$

(II) $[\overline{\xi}(\alpha), \underline{\xi}(\alpha)] - [\overline{\gamma}(\alpha), \underline{\gamma}(\alpha)] = [\min\{\lim_{n\to\infty} m_1 - \lim_{n\to\infty} m_2, \lim_{n\to 1} m_1 - \lim_{n\to 1} m_2\}, \max\{\lim_{n\to\infty} m_1 - \lim_{n\to\infty} m_2, \lim_{n\to 1} m_1 - \lim_{n\to 1} m_2\}]$

(III) $[\overline{\xi}(\alpha), \underline{\xi}(\alpha)] \times [\overline{\gamma}(\alpha), \underline{\gamma}(\alpha)] = [\min\{\lim_{n\to\infty} m_1 \times \lim_{n\to\infty} m_2, \lim_{n\to 1} m_1 \times \lim_{n\to 1} m_2\}, \max\{\lim_{n\to\infty} m_1 \times \lim_{n\to\infty} m_2, \lim_{n\to 1} m_1 \times \lim_{n\to 1} m_2\}]$

(IV)  $[\overline{\xi}(\alpha),\underline{\xi}(\alpha)] \div [\overline{\gamma}(\alpha),\underline{\gamma}(\alpha)] = [\min\{\lim_{n\to\infty} m_1 \div \lim_{n\to\infty} m_2, \lim_{n\to 1} m_1 \div \lim_{n\to 1} m_2\}, \max\{\lim_{n\to\infty} m_1 \div \lim_{n\to\infty} m_2, \lim_{n\to 1} m_1 \div \lim_{n\to 1} m_2\}]$

In Fig. 2, schematic diagram of FFEM is presented. It has three layers viz. input, hidden layer, and output. In the input layer, the involved system parameters are taken as TFN then that numbers are executed in the hidden layer by using FEM algorithm. Finally, the solutions are obtained as output.

In Fig. 3, the FEM algorithm is described which is used in the hidden layer given in Fig. 2.

We take the convection diffusion equation with the Galerkin residue approach and finite element procedure to solve the governing differential equation. The final form of the residue in matrix form can be written as

$$[K]\{T\} = \{f\}, \tag{9}$$

were,

$$[K] = \int_\Omega [B]^T[D][B]d\Omega + \int_\Gamma h[N]^T[N]d\Gamma \tag{10}$$

$$\{f\} = \int_\Omega G[N]^T d\Omega - \int_\Gamma q[N]^T d\Gamma + \int_\Gamma hT_\infty[N]^T d\Gamma. \tag{11}$$

Here, $K$ is the global stiffness matrix, $T$ presents the temperature, $f$ is the force vector, and $N$ is defined as the shape functions.

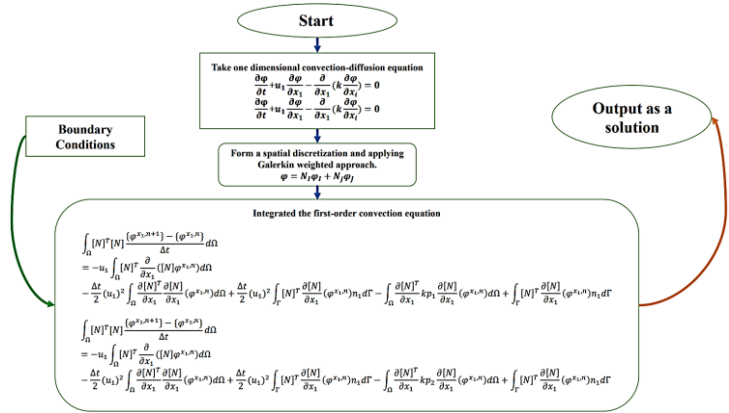

Fig. 3. Flow chart for Finite Element Method (FEM) to solve transient convection-diffusion problem

## 3. Transient convection-diffusion heat transfer on a square plate

The above discussed algorithm is used to solve convection diffusion problem for a rectangular plate. We have considered a two-dimensional rectangular plat having dimension $20cm \times 10cm$. Where the left wall of the plate specified heat flux and the right wall maintained constant temperature. The top wall is convective heat transfer and down wall is adiabatic. The parameters are, $h$- the convection heat transfer, $q$- heat input rate and ambient temperature $T_a$ are considered 10 percentage of error uncertainty. Initially, the temperature is investigated by taking crisp value by using traditional finite element method. Next the parameters are taken as uncertain in left and right side of fuzzy number, then the temperature are obtained by taking Fuzzy Finite Element Method. Here the plate is discretised into 50 triangular elements. After applying boundary conditions, we got the left and right elemental temperatures. The same is presented in Figs. 4 and 5. Respectively. Tables 1 provide the different uncertain parametric values.

Table 1. Crisp and fuzzy values of involved parameters

| Parameters | Crisp value | Uncertain parametric values |
|---|---|---|
| $h$ | $1.2 W/cm^2 K$ | $1.2 \pm 5\% \; W/cm^2 K$ |
| $q$ | $2 W/cm^2$ | $2 \pm 5\% \; W/cm^2$ |

Now we discuss the variation of uncertain temperature which is obtained by taking different combination of fuzzy parameters. First, we choose the convection heat transfer as TFN and then the problem is solved using mentioned FFEM. The obtained results are depicted in Fig. 4. Similarly, taking only heat input rate as TFN the obtained temperatures are shown in Fig. 5.

To check the sensitiveness of the uncertain parameters, the uncertain width is considered for each case (only $h$ is fuzzy and only $q$ is fuzzy) then the obtained widths are compared and reported in Figs. 6 and 7. It is found that the average width is 0.374762632 when only $h$ is fuzzy, whereas for only $q$ is fuzzy, the average width is 0.469378335. Similarly, the variances are 0.045947279 and 0.193681158 for only $h$ is fuzzy and only $q$ is fuzzy respectively. Further it is observed that the case when only $q$ is fuzzy, system is more sensitive. Therefore, small change in $q$ results more uncertain solutions.

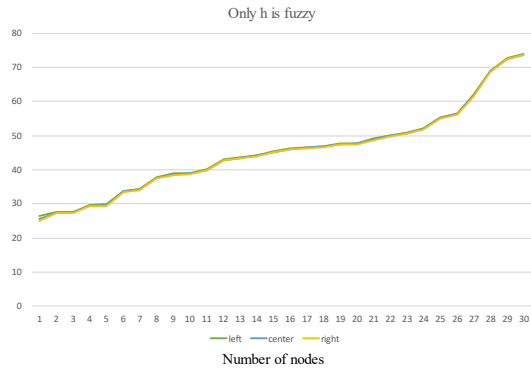

Fig. 4. Uncertain distribution of temperatures when convective heat transfer ($h$) is fuzzy

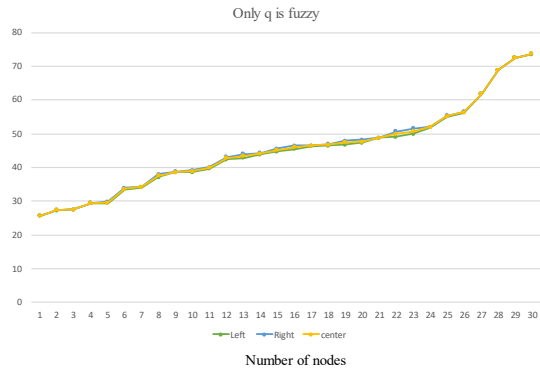

Fig. 5. Uncertain distribution of temperatures when heat input rate ($q$) is fuzzy

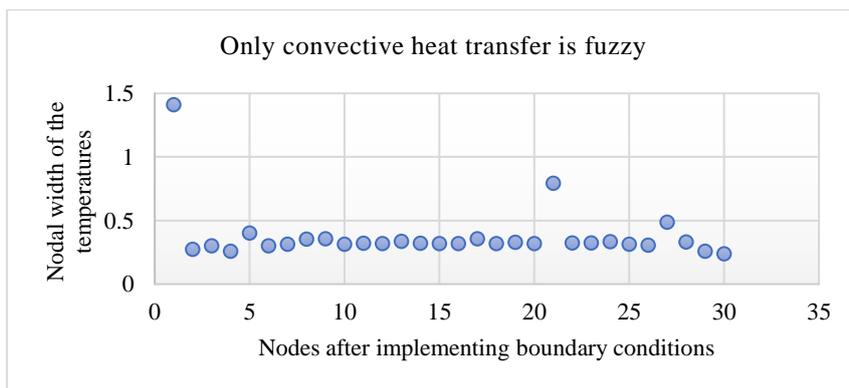

Fig. 6. Variation of uncertain temperatures when convective heat transfer is fuzzy

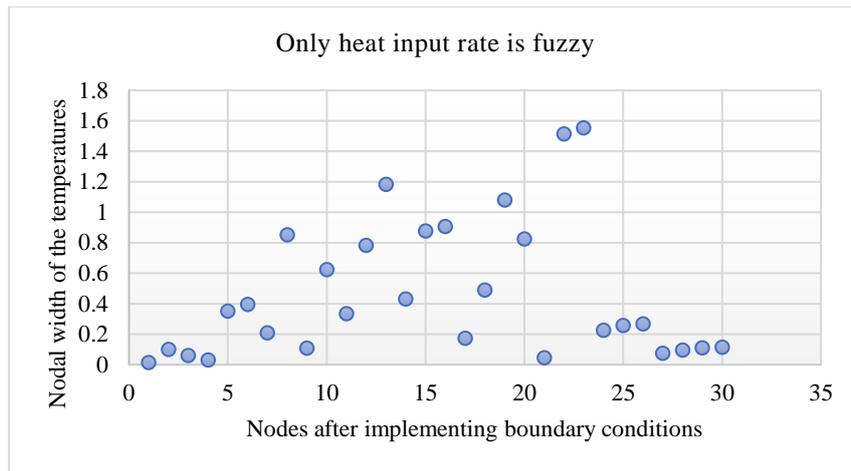

Fig. 7. Variation of uncertain temperatures when heat input rate is fuzzy

## 4. Conclusion

This paper presents a theoretical approach to investigated the different combination of uncertain parameters in the convection-diffusion equation and the sensitiveness of the parameters. A fuzzy finite element model is modelled for transient convection-diffusion problems. The parametric representation of fuzzy numbers to crisp form is adopted to compute the uncertainties that ease the process. The hidden layer of the overall process is performed through FEM. Then the said algorithm is used for a case study and for the mentioned problem it is noticed that $q$ is more sensitive than $h$ to the system.